# Data-driven optimization of reliability using buffered failure probability

**Ji-Eun Byun [1*] and Johannes O. Royset [2]**

[1] Department of Civil, Environmental and Geomatic Engineering, University College London, London, United Kingdom; j.byun@ucl.ac.uk

[2] Operations Research Department, Naval Postgraduate School, California, United States; joroyset@nps.edu

[*] Correspondence: j.byun@ucl.ac.uk; Tel.: +44-20-3108-1501

**Abstract:** Design and operation of complex engineering systems rely on reliability optimization. Such optimization requires us to account for uncertainties expressed in terms of complicated, high-dimensional probability distributions, for which only samples or data might be available. However, using data or samples often degrades the computational efficiency, particularly as the conventional failure probability is estimated using the indicator function whose gradient is not defined at zero. To address this issue, by leveraging the *buffered failure probability*, the paper develops the *buffered optimization and reliability method* (BORM) for efficient, data-driven optimization of reliability. The proposed formulations, algorithms, and strategies greatly improve the computational efficiency of the optimization and thereby address the needs of high-dimensional and nonlinear problems. In addition, an analytical formula is developed to estimate the reliability sensitivity, a subject fraught with difficulty when using the conventional failure probability. The buffered failure probability is thoroughly investigated in the context of many different distributions, leading to a novel measure of tail-heaviness called the *buffered tail index.* The efficiency and accuracy of the proposed optimization methodology are demonstrated by three numerical examples, which underline the unique advantages of the buffered failure probability for data-driven reliability analysis.

**Keywords:** Reliability optimization; data-driven optimization; buffered failure probability; superquantile; tail index; reliability sensitivity

## 1 Introduction

In order to secure the resilience of real-world engineering systems (e.g., structural systems, infrastructures, and mechanical systems), risk-informed decisions should be made by correctly accounting for uncertainties arising from natural and humanmade hazards, external loads, material properties, and various other sources (Fenton and Neil 2018). Such decisions can be supported by the solution of optimization problems that identify the least expensive decisions with acceptable failure probability. These optimization problems tend to be computationally challenging when dealing with nonlinear limit-state functions as well as many random variables and decision variables. If the corresponding probability distributions are unknown or rather complex, then it becomes necessary to leverage data and/or samples generated from models of uncertainty. However, data and samples tend to cause challenges for optimization problems involving failure probability constraints. In the face of the recent advancement of data technology, it is essential to be able to carry out reliability optimization using data or samples only.

The conventional probability is expressed as a weighted sum of indicator functions in settings with only data or samples, and this causes significant challenges for optimization algorithms. Specifically, the indicator function "counts" the number of failure events, but in the process, introduces expressions that are not differentiable; the indicator function lacks a derivative at zero. All subsequent inference tasks involving derivatives and gradients, as in sensitivity analysis and optimization, then become challenging and/or poorly defined.

Reliability sensitivity is typically defined as the derivative of a failure probability with respect to the parameter or variable of interest. It is useful for quantifying the impacts of variables on the failure probability. To avoid expressions involving the indicator functions, studies employ costly simulation methods that monitor the boundary of limit-state functions such as directional sampling (Royset and Polak 2007; Ackooij and Henrion 2017; Valdebenito et al. 2018). Another approach is to faciliate computations by approximating problems with differentiable functions such as linearizations of response surfaces (Melchers and Ahammed 2004) and smoothing of indicator functions (Papaioannou et al. 2018; Kannan and Luedtke 2020).

Commonly used reliability methods such as FORM and SORM require us to transform the random variables into the space of standard normal random variables (Der Kiureghian 2005), which may incur errors or even become inapplicable for distributions with highly dependent variables or non-conventional formulas (Liu and Der Kiureghian 1986). Accordingly, for data-driven reliability analysis, such approximation is likely to become impossible when the given data are high-dimensional or deviate significantly from any of the common distributions. Another disadvantage of such transformation is that the inference over random variables becomes complicated as they reside in a different space from the one used to define limit-state functions, parameters, and decision variables. For example, reliability sensitivity with respect to a random variable requires additional treatment to take into account the relation between the two spaces (Papaioannou et al. 2018), while the computation becomes even more complicated when random variables and other variables need to be considered together.

To address such difficulties in formulations and computations associated with the conventional failure probability, an alternative measure of reliability has been proposed, namely the *buffered failure probability* (Rockafellar and Royset 2010). While the buffered failure probability has received a limited attention in reliability engineering, its basic ideas have been actively applied in finance, economics, and operations research through the closely related concept *conditional value-at-risk (CVaR)* (Rockafellar and Uryasev 2002), which is also called *average value-at-risk, expected shortfall,* and *superquantile*. In contrast to the conventional failure probability that defines failure as an event producing a positive value of the limit-state function, the buffered failure probability relies on a threshold above which the limit-state function has an average value of zero. The primary advantage of this alternative idea is that it does not require the indicator function and thereby facilitates inference tasks that remain challenging in the context of the conventional failure probability. This advantage has been demonstrated with various inference tasks such as portfolio management (Rockafellar and Uryasev 2002), binary classifications (Norton and Uryasev 2019), design optimization of structural systems (Rockafellar and Royset 2010), and evaluation of derivatives (Zhang et al. 2019; He 2020). Still, it is noted that this alternative definition does not nullify the established standards for the conventional probability as the two probabilities in general show positive correlation (i.e., a larger value of either probability implies a larger value of the other one).

The secondary advantage of the buffered failure probability is that it reflects not only the probability of exceeding a threshold, but also the tail's heaviness. In other words, given a pair of distributions with the same conventional failure probability, the buffered probability is larger for the distribution with the heavier tail. The utility of accounting for tail behaviors has been widely recognized and used to evaluate tail-related risk, particularly by employing the concept of CVaR (Bekiros et al. 2019; Echaust and Just 2020). The tail of a distribution represents high-consequence failure events and thus needs to be considered in a decision-making context. Heaviness of tails can also be measured using parametric approaches (De Haan and Ferreira 2007; Németh and Zempléni 2020) or samples (Qi 2010; Xiong and Peng 2020). While these existing measures provide useful insights, they may show a limited applicability by being developed based on a certain class of distributions or lacking extensive investigations over general problems.

CVaR has been recognized to possess certain advantages and potentials in the areas of finance and economics (Dixit and Tiwari 2020), industrial engineering (Rezaei et al. 2020; Zhu et al. 2020), energy engineering (Feng et al. 2020), machine learning (Soma and Yoshida 2020), and design optimization (Li et al. 2021). However, the buffered failure probability has been studied only to a limited extent in reliability engineering (Rockafellar and Royset 2010; Zrazhevsky et al. 2020; Chaudhuri et al. in review). Accordingly, to demonstrate the relevance of this alternative reliability measure for engineering systems, this paper aims to investigate the associated practical issues and methods for addressing reliability problems in a data-driven setting. This paper presents a novel formula for the buffered failure probability in the data-driven setting, constructs a new expression for its derivative, investigates the relation to the conventional failure probability, and quantifies distributions with heavy tails. Moreover, we develop a novel algorithm for optimization of the buffered failure probability based on an adaptive sample size and an active-set strategy. We term the general approach of the developed formulas and algorithms the *buffered optimization and reliability method* (BORM).

This paper is organized as follows. Section 2 introduces the buffered failure probability and derives a formula to compute the probability value using discrete outcomes, samples, or data. This section also compares computational advantages and practical issues with the conventional probability. Then, to better understand this alternative reliability measure, Section 3 derives and analyzes the ratio between the two probabilities using common distributions and proposes an index that quantifies the heaviness of a distribution's tail. Reference values and thresholds of the index are given as well so that its implementation can be facilitated. Section 4 develops BORM; it illustrates the formulations of data-driven reliability optimization using buffered failure probability and proposes practical algorithms and strategies by which the computational cost of the optimization can become affordable. In addition, Section 5 develops an approach for estimating the reliability sensitivity of buffered failure probability using samples or data. The applicability and accuracy of BORM are demonstrated by numerical examples in Section 6, whose computational efficiency is underlined by being solved with a personal desktop and general-purpose optimization solvers provided by Matlab®. Specifically, two benchmark examples demonstrate that compared to the existing methods, BORM obtains better solutions through its mathematically better behaved properties and maintains computational cost at practical levels. A truss bridge example demonstrates the unique advantage of the proposed optimization methodology for data-driven decision-making. The supporting source code and data are available for download at (URL which will be provided here once the paper is accepted.)

While BORM can be coupled with any advanced sampling techniques, which would reduce the required number of samples, the following discussion and numerical examples focus on Monte Carlo Simulation (MCS) for two reasons: (1) Since MCS is straightforward to implement, the proposed algorithms can remain accessible and practical. (2) As MCS is representative of naturally acquired data, e.g., weigh-in-motion (WIM) data or wind load data, the proposed algorithms are suitable for data-driven optimization, which is particularly useful in view of the recent advancement of data technology.

Throughout the paper, the conventional and buffered failure probabilities are denoted by $p_f$ and $\bar{p}_f$, respectively. The letters $V$ and $x$ stand for random and decision variables, which are bolded when referring to a vector of variables, i.e., $\boldsymbol{V}$ and $\boldsymbol{x}$. The two types of variables together determine the value of the limit-state function $g(\boldsymbol{x}, \boldsymbol{V})$, while the conventional failure event is defined by $g(\boldsymbol{x}, \boldsymbol{V}) > 0$. Although it is also common to define the failure event as $g(\boldsymbol{x}, \boldsymbol{V}) < 0$, the aforementioned definition is adopted so that the following developments are compatible with the preceding works (Rockafellar and Royset 2010; Mafusalov and Uryasev 2018; Norton et al. 2019; Royset and Wets 2021). A conversion is easily accomplished by reversing the sign.

## 2 Buffered failure probability

This section recalls the definitions of the buffered failure probability and discusses its computation based on samples and data. This leads to an expression in Section 2.1. The discussion is limited to the materials that are directly related to the scope of the paper. Additional theoretical properties of the buffered failure probability can be found in Rockafellar and Royset (2010), Mafusalov and Uryasev (2018), and Section 3.E of Royset and Wets (2021).

### 2.1 Background

Given a random variable $Y$ with cumulative distribution function (CDF) $F_Y(\cdot)$, the conventional failure probability is defined as

$$p_f = P[Y > 0] = 1 - F_Y(0). \tag{1}$$

If $F_Y(y)$ is strictly increasing, then the *α-quantile* for $\alpha \in (0,1)$, denoted by $q_\alpha$, is given as[1]

$$q_\alpha = F_Y^{-1}(\alpha). \tag{2}$$

Figure 1(a) illustrates the situation for $\alpha = 0.84$ and a normal distribution with mean −1 and standard deviation 1. Then, one obtains the equivalent expression

$$p_f = 1 - \alpha_0, \tag{3}$$

where $\alpha_0$ is the probability that makes the quantile equal to zero, i.e.,

$$q_{\alpha_0} = 0. \tag{4}$$

In contrast, the buffered failure probability is obtained by replacing the α-quantile in the formula for $p_f$ by the mathematically better behaving *α-superquantile*. The α-superquantile of $Y$, denoted by $\bar{q}_\alpha$, is a modification of the α-quantile of $Y$ by accounting for the average value of the outcomes beyond the α-quantile. Specifically, it is given as[2]

$$\bar{q}_\alpha = q_\alpha + \frac{1}{1-\alpha} \mathbb{E}[\max\{Y - q_\alpha, 0\}]. \tag{5}$$

---

[1] For general $F_Y(y)$, the α-quantile equals the smallest value $y$ with $F_Y(y) \geq \alpha$.

[2] In the finance and operations research literature, a superquantile is called CVaR (Rockafellar and Uryasev 2002). Here, the terminology follows Rockafellar and Royset (2010).

Then, the buffered failure probability is defined as

$$\bar{p}_f = 1 - \bar{\alpha}_0, \tag{6}$$

where $\bar{\alpha}_0$ is the probability that makes the superquantile equal to zero, i.e.,

$$\bar{q}_{\bar{\alpha}_0} = 0. \tag{7}$$

This definition of $\bar{p}_f$ applies as long as $\mathbb{E}[Y] < 0$ and $p_f > 0$, which is hardly a limitation for practical problems. In the rare situation with $p_f = 0$, one defines $\bar{p}_f = 0$. In all other cases, one sets $\bar{p}_f = 1$. If $Y$ is continuously distributed, it always holds that $F_Y(q_\alpha) = \alpha$, and the expressions for $\bar{p}_f$ and $\bar{q}_\alpha$ are simplified to

$$\bar{q}_\alpha = \mathbb{E}[Y|Y \geq q_\alpha] \tag{8}$$

and

$$\begin{aligned}\bar{p}_f &= 1 - \bar{\alpha}_0 \\ &= 1 - F_Y(q_{\bar{\alpha}_0}) \\ &= P[Y > q_{\bar{\alpha}_0}].\end{aligned} \tag{9}$$

For the normal random variable with mean $\mu = -1$ and standard deviation $\sigma = 1$ in Figure 1, the superquantile can be evaluated by (8). In this case,

$$\mathbb{E}[Y|Y \geq y] = \mu + \frac{\sigma \cdot \phi\left(\frac{y-\mu}{\sigma}\right)}{1 - \Phi\left(\frac{y-\mu}{\sigma}\right)}, \tag{10}$$

where $\phi(\cdot)$ and $\Phi(\cdot)$ are the probability density function (PDF) and CDF of the standard normal distribution; see Figure 1(b). From the figure, it is concluded that $y = -0.71$ produces $\mathbb{E}[Y|Y \geq y] = 0$. In turn, $-0.71$ is the $\alpha$-quantile of $Y$ when $\alpha = 0.61$. Thus, by (8) and (9), $\bar{\alpha}_0 = 0.61$ and $\bar{p}_f = 0.39$; Figure 1(c) visualizes these probabilities using the PDF of $Y$.

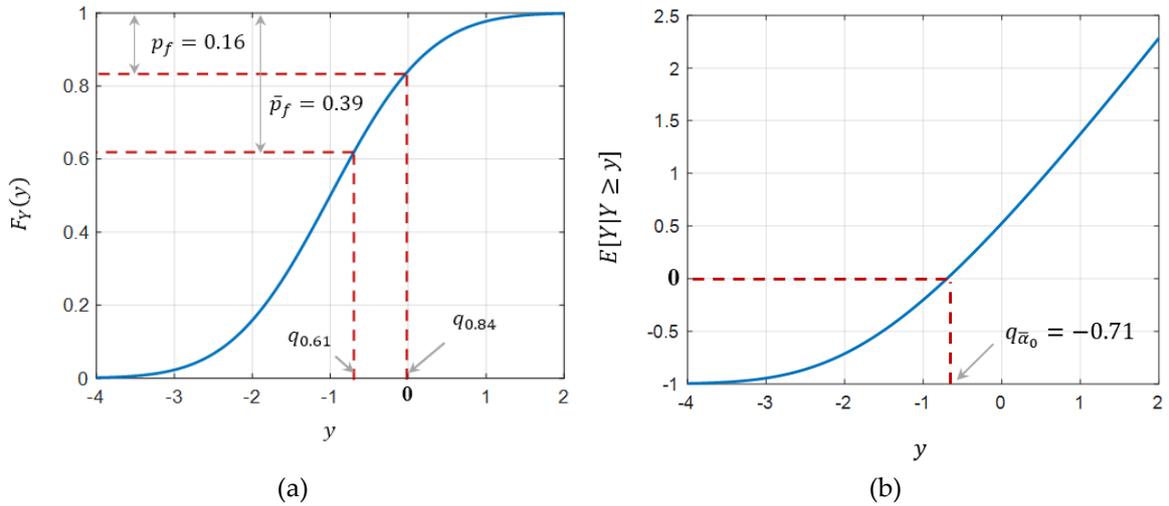

(a)  (b)

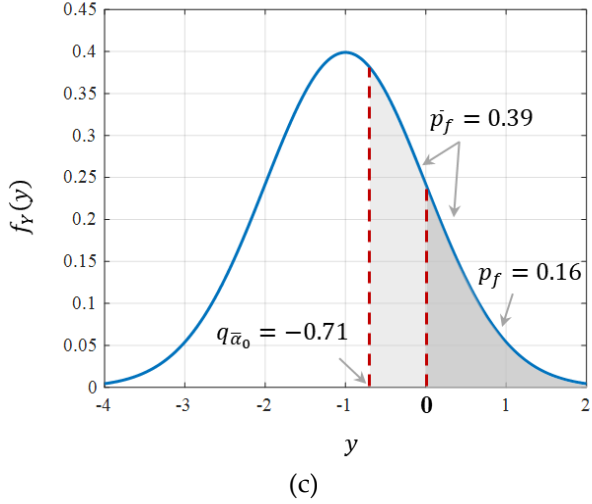

(c)

Figure 1. Buffered failure probability of the normal distribution with mean −1 and standard deviation 1: (a) CDF and failure probabilities, (b) conditional expectation, and (c) PDF and failure probabilities.

## 2.2 Closed-form formula using samples or data

Leveraging Mafusalov and Uryasev (2018), this section derives a formula for the buffered failure probability $\bar{p}_f$ in the case of a random vector $\boldsymbol{V}$ with outcomes $\boldsymbol{v}_1, \cdots, \boldsymbol{v}_N$ occurring with probabilities $p_1, \cdots, p_N$. If the outcomes are generated by sampling, the probabilities correspond to sample weights; for example, if they are generated by Monte Carlo Simulation (MCS), $p_n = 1/N$ for all $n = 1, \cdots, N$. For a fixed design vector $\boldsymbol{x}$, the limit-state function has the outcomes $y_n = g(\boldsymbol{x}, \boldsymbol{v}_n)$, $n = 1, \cdots, N$. Without loss of generality, it is assumed that $y_n < y_{n'}$ if $n < n'$ as one can always reorder the outcomes $\boldsymbol{v}_1, \cdots, \boldsymbol{v}_N$. Moreover, if two outcomes $y_n$ and $y_{n'}$ are identical, they can be simply treated as a single outcome with probability $p_n + p_{n'}$. Then,

$$\bar{p}_f = \sum_{n=n^*+1}^{N} \left( \frac{p_n (y_n - y_{n^*})}{-y_{n^*}} \right), \tag{11}$$

where $n^*$ is an integer in $\{1, \cdots, N\}$ that satisfies

$$\sum_{n=n^*}^{N} p_n y_n < 0 < \sum_{n=n^*+1}^{N} p_n y_n. \tag{12}$$

This expression for the buffered failure probability is novel and follows by simplifying several aspects of a formula furnished by Corollary 2.4 in Mafusalov and Uryasev (2018). It is straightforward to use: Start by checking (12) with $n^* = N - 1$. If the inequalities do not hold, then check $n^* = N - 2$ and so forth. Thus, the computational complexity of computing $\bar{p}_f$ is $O(N \log N)$, with the sorting of the values $y_1, \cdots, y_N$ being the bottleneck; see for example, Cormen et al. (2009). While the formula addresses essentially all practical situations, there are pathological cases when it does not apply. These cases are characterized by having no $n^*$ satisfying (12). For example, this would be the case if all $y_n < 0$. One can then typically determine the buffered failure probability directly from its definition.

The above expression for the buffered failure probability is for limit-state functions with a finite number of outcomes. In the case of a continuously distributed or high-dimensional discrete random vector $\boldsymbol{V}$, $g(\boldsymbol{x}, \boldsymbol{V})$ would have an insurmountable number of outcomes. Then, one can still bring in the formula as an estimator of the buffered failure probability.

Simply sample $V$ to produce the outcomes $v_1, \cdots, v_N$ and set $p_n$ as the corresponding sample weights. Compute $y_n = g(x, v_n)$ and sort the values. As above, remove any duplicates and update the probabilities, i.e., if $y_n = y_{n'}$, then use $p_n + p_{n'}$ for this outcome. Thus, we obtain the estimate

$$\bar{p}_f \approx \hat{\bar{p}}_f = \sum_{n=n^*+1}^{N} \left( \frac{p_n(y_n - y_{n^*})}{-y_{n^*}} \right). \tag{13}$$

## 2.3 Computational advantages for data-driven reliability optimization

In the data-driven setting, we are faced with a finite number of outcomes $y_1, \cdots, y_N$ as in Section 2.2. Then, the buffered failure probability $\bar{p}_f$ improves computational efficiency for two reasons. First, it is noted that the conventional failure probability $p_f$ is evaluated as

$$p_f = \sum_{n=1}^{N} p_n \cdot \mathbb{I}(y_n > 0), \tag{14}$$

where $\mathbb{I}(\cdot)$ is the indicator function that takes value 1 if the given statement is true and 0 otherwise. Since the derivative of the indicator function is not defined at zero, the computation of gradients with respect to the parameters in the limit-state function becomes ill-defined or ineffective. This has ramification for various inference tasks including optimization and sensitivity estimation. In contrast, as observed in (11) and (13), the indicator function is not used for computing $\bar{p}_f$. Accordingly, $\bar{p}_f$ is mathematically better behaved and leads to improved computational efficiency as demonstrated in Sections 4 and 5.

Second, while the difficulties arising from the indicator function can be circumvented by reliability methods such as the first- and second-order reliability method (FORM and SORM), those methods cause other challenges related to differentiability (Rockafellar and Royset 2010) and also require us to transform the random variables to the standard normal space (Der Kiureghian 2005). Although they have been effectively implemented, transformed domains might not preserve some characteristics of the original space; and furthermore, such transformation becomes inapplicable when dealing with highly correlated random variables or nonconventional distributions (Liu and Der Kiureghian 1986). In contrast, transformations are completely avoided by utilizing $\bar{p}_f$, whereby a wider class of distributions can be addressed. In the case of data-driven optimization, $\bar{p}_f$ exempts us from fitting the data to preselected probability distributions, which is particularly challenging when dealing with high-dimensional datasets. Another advantage is that since random variables remain in the original space, the inferences do not need to distinguish between random variables and other variables (e.g., decision variables and the parameters of limit-state functions), as discussed in detail in Section 5.

## 2.4 Invariance issue

One of the major differences between $p_f$ and $\bar{p}_f$ is the threshold used to define failure. In the former case, the threshold is fixed at zero. The latter probability can be interpreted as having a "flexible" threshold $q_{\bar{\alpha}_0}$ that is adjusted based on the tail of the underlying probability distribution; see (9). As a result, while $p_f$ remains the same with varying formulas of a limit-state function as long as the failure domain remains the same, $\bar{p}_f$ could change. For example, given random variables $S$ and $R$, the two limit-state functions $g_1(s, r) = s - r$ and $g_2(s, r) = s/r - 1$ lead to the same failure domain and $p_f$, i.e., $\{s, r | s - r > 0\} =$ $\{s, r | s/r - 1 > 0\}$, but typically to different $\bar{p}_f$ because of the potentially different tails of the distributions of $S - R$ as compared to $S/R - 1$.

The lack of invariance is not necessarily a weakness for two reasons. First, this implies that in order to use $\bar{p}_f$, engineering judgements should be made to correctly define the measure of interest, which is typically a simple task. For example, if $S$ and $R$ stand for the external loads and structural strength, respectively, $g_1$ represents absolute load exceedance, while $g_2$ implies relative load exceedance. Then, one should be able to make a proper choice between the two measures based on the consequence of interest. Second, the invariance property has been exploited to facilitate the computation of $p_f$, e.g., when $S$ and $R$ follow the lognormal distribution, $g_1$ can be transformed to $g_2$ so that $p_f$ can be evaluated analytically. However, such advantage is not relevant to the current study where the focus is on reliability optimization using samples or data, rather than the analytical calculation of reliability.

## 3 Buffered tail index and target probability calibration

While there are widely accepted norms for choosing a permissible level of the conventional failure probability $p_f$ in terms of the reliability index and return periods, such standards are not available for the buffered failure probability $\bar{p}_f$ owing to its relatively short history of usage. To address this issue, we investigate the relationship between $p_f$ and $\bar{p}_f$ using common distributions and thereby facilitates the translation of permissible levels of $p_f$ into permissible levels of $\bar{p}_f$. From the investigations, we propose to use the ratio of the two probabilities as a measure of the heaviness of a distribution's upper tail, i.e.,

$$\tau = \bar{p}_f / p_f, \qquad (15)$$

which is coined the *buffered tail index*. This measure provides useful insights for practical decision-making as upper tails are representative of the events of primary concern; those with low probability and high consequence. Moreover, the buffered tail index leads to a useful rule for translating a desired level of $p_f$ into a bound on $\bar{p}_f$ in the context of reliability optimization.

In order to understand how the shape of a distribution determines $\tau$, we investigate common distributions and compare them with two reference distributions, the normal and exponential distributions (Explicit formulas for $\tau$ are given in Appendix A). Figures 2 and 3 illustrate distributions with light and heavy upper tails, respectively, where the parameters of the distributions are set to satisfy $p_f = 0.1$. Figure 2 investigates the lognormal distribution with $s = 0.125$, Weibull distribution with $k = 1.5$, and GEV with $\xi = 0$. In Figure 2(a), the $\tau$ values of these distributions lie between those of the two reference distributions; and Figures 2(b) and (c) illustrate their PDFs and upper tails. In Figure 2(c), it is noted that the upper tails of these distributions lie in the vicinity of the tails of the two reference distributions. Figure 3 illustrates similar investigations for the lognormal distribution with $s = 1$, Weibull distribution with $k = 0.5$, and GEV with $\xi = 0.5$. Figure 3(a) shows that these distributions have greater $\tau$ values than the reference distributions since they have heavier upper tails in Figures 3(b) and (c).

Figures 2 and 3 confirm that the buffered tail index $\tau$ quantifies the level of right-tail heaviness for a distribution. To determine the boundary value of $\tau$ between light and heavy tails, it is noted that heavy-tailed distributions are defined as those with tails that are not exponentially bounded, i.e., they have a heavier tail than the exponential distribution (Embrechts et al. 2013). Therefore, it is proposed to set the threshold as the $\tau$ value of the exponential distribution, $e \approx 2.72$, i.e., distributions with $\tau > 2.72$ are classified as heavy-tailed distributions. For example, $\tau = 4$ indicates a distribution with a significant possibility of high-consequence failure events that warrant further investigation.

(a)
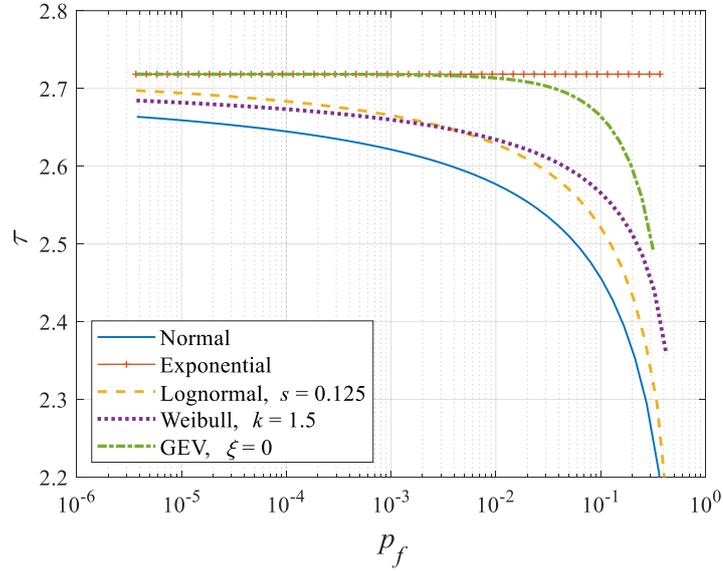

(b)
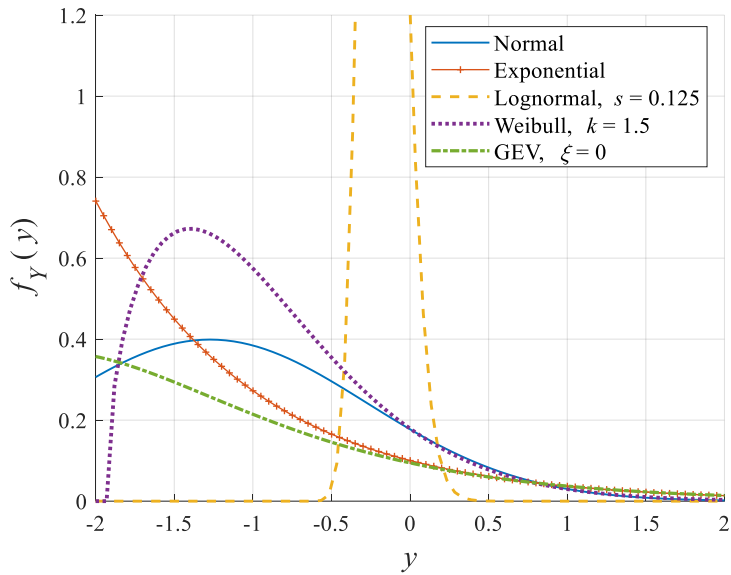

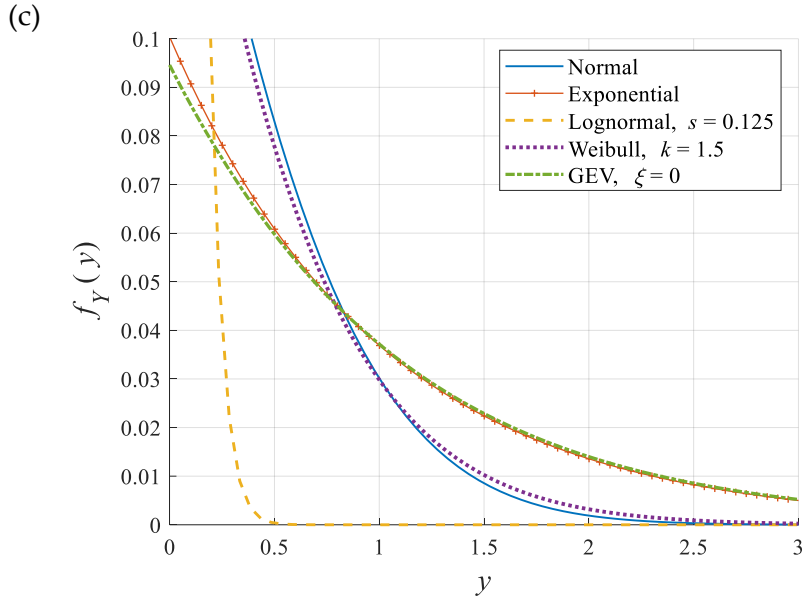

Figure 2. Distributions with light upper tails: lognormal with $s = 0.125$, Weibull with $k = 1.5$, and GEV with $\xi = 0$. (a) Buffered tail index, (b) PDFs with $p_f = 0.1$, and (c) upper tails of the PDFs.

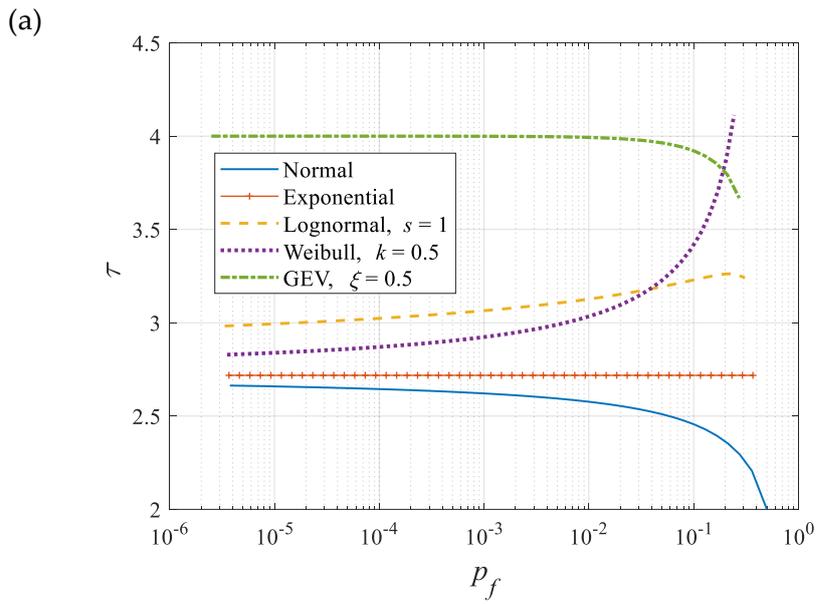

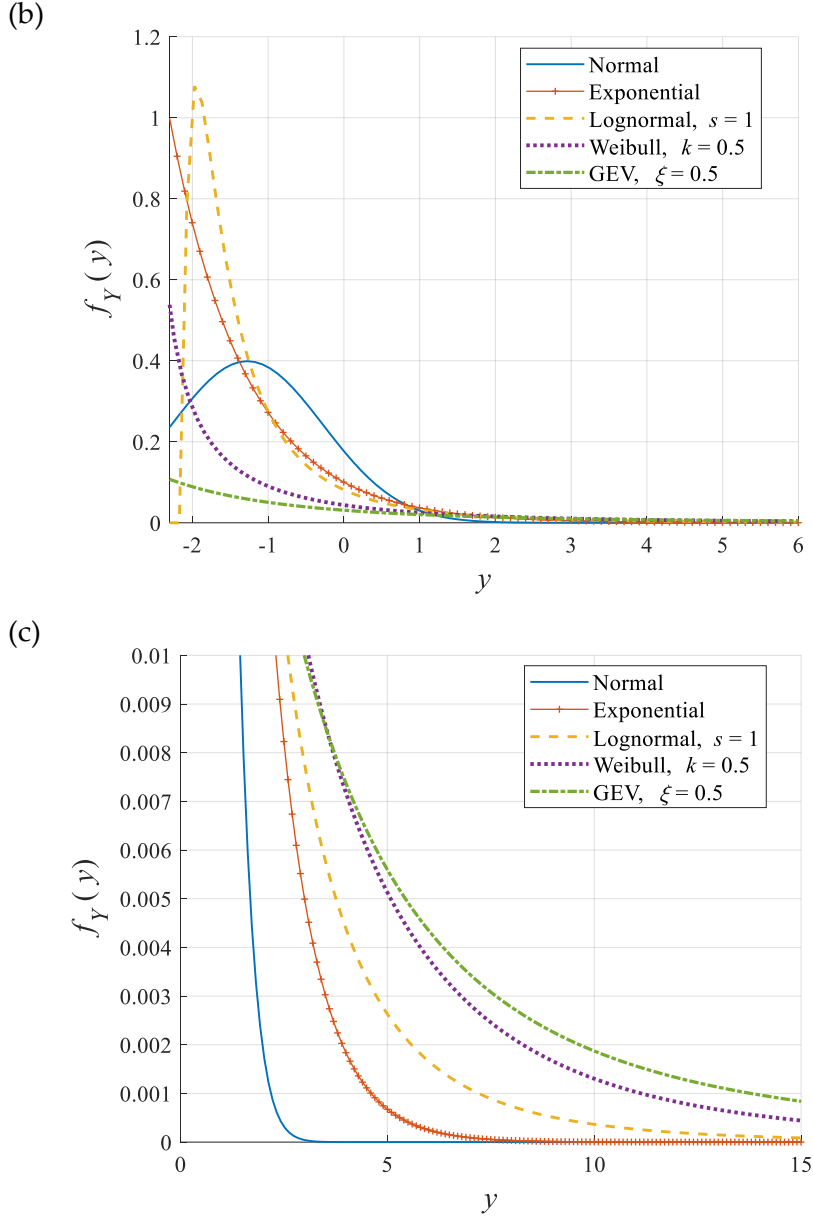

Figure 3. Distributions with heavy upper tails: lognormal with $s = 1$, Weibull with $k = 0.5$, and GEV with $\xi = 0.5$. (a) Buffered tail index, (b) PDFs with $p_f = 0.1$, and (c) upper tails of the PDFs.

In order to perform reliability optimization, a target buffered failure probability $\bar{p}_f^t$ needs to be selected. To inform the selection, one can leverage the buffered tail index $\tau$ and a target $p_f^t$ for the conventional probability. The latter is typically based on a reliability index (e.g. 2) or a return period (e.g. 100 years) (Kim and Song 2021). To this end, it is proposed to use the $\tau$ values of the normal distribution as the criterion for two reasons: (1) as observed in Figures 3 and 4, the normal distribution has the lowest $\tau$ and thereby, leads to the most conservative estimation of $\bar{p}_f^t$ given the same level of $p_f^t$, and (2) it is the most broadly used distribution, especially when the true distribution is unknown. Figure 4 illustrates the $\tau$ values of the normal distribution, $\tau_{\text{normal}}$. In the figure, the analytical derivation (solid blue line) agrees with the results computed using $10^7$ samples (dashed orange line). It is noted that $\tau_{\text{normal}}$ depends only on $p_f$ and not on the distribution parameters ($\mu$ and $\sigma$).

Based on these observations, the reference value of $\tau$ can be set to the fitted approximation for $p_f \in [10^{-6}, 0.5]$ in Figure 4 (dotted yellow line), i.e.,

$$\tau^* = \begin{cases} \dfrac{2.61 - 2.68}{\ln(0.01) - \ln(10^{-6})}\left(\ln p_f - \ln(10^{-6})\right) + 2.68, & 10^{-6} \leq p_f < 0.01 \\ \dfrac{2.4 - 2.61}{\ln(0.3) - \ln(0.01)}\left(\ln p_f - \ln(0.01)\right) + 2.61, & 0.01 \leq p_f < 0.3 \\ \dfrac{2 - 2.4}{\ln(0.5) - \ln(0.3)}\left(\ln p_f - \ln(0.3)\right) + 2.4, & 0.3 \leq p_f \leq 0.5 \end{cases} \qquad (16)$$

This expression can be used to determine the target buffered failure probability as $\bar{p}_f^t = \tau^* \cdot p_f^t$ in reliability optimization. In Section 6, the proposed conversion is found to be useful. Still, some adjustments can be expected in cases with especially high or low $\tau$.

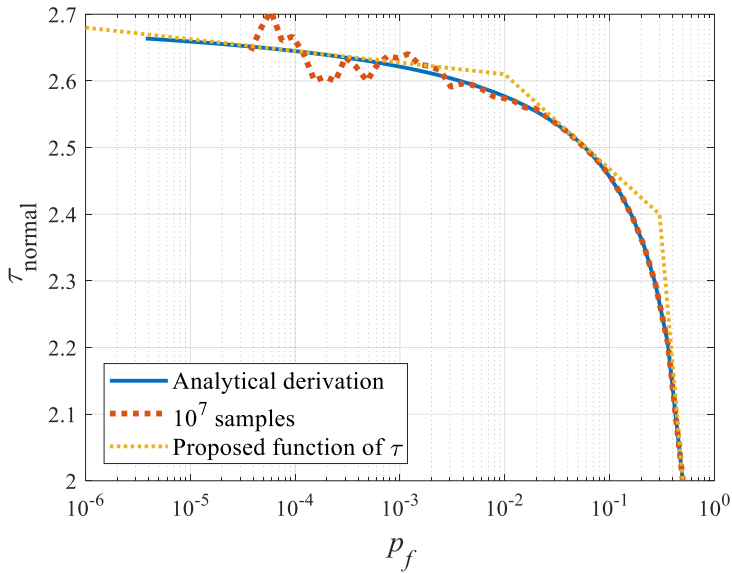

Figure 4. Buffered tail index of the normal distribution and the proposed reference formula of $\tau$.

## 4 Data-driven optimization using buffered failure probability

In this section, practical formulations and algorithms are proposed for reliability optimization using samples or data. We refer to the overall approach as the *buffered optimization and reliability method* (BORM). Section 4.1 develops the basic formulations of the optimization problem. Then, to facilitate implementation, the following sections address practical issues: Section 4.2 develops an algorithm for optimization to reduce the computational cost via an active-set strategy; Section 4.3 discusses handling multiple limit-state functions and finding a feasible solution so that an initial solution can be provided before starting the optimization. Finally, Section 4.4 summarizes the proposed optimization procedure.

### 4.1 Formulations of optimization problem

A reliability optimization problem aims to minimize a cost function while ensuring a failure probability below a target value, i.e.,

$$\min_{x} c(x)$$
$$\text{subject to } p_f(x) \leq p_f^t \tag{17}$$
$$x \in X,$$

where $c(\cdot)$, $p_f^t$, and $X$ denote the cost function, target $p_f$, and a set that represents additional constraints imposed on the decision variables, respectively. In the optimization problem, the conventional failure probability is evaluated as a function of $x$, i.e.,

$$p_f(x) = P[g(x, V) > 0]. \tag{18}$$

Given discrete outcomes $v_n$ and corresponding probabilities $p_n$, $n = 1, \cdots, N$, that are specified by the distribution of random variables $V$, the optimization problem stated in (17) becomes

$$\min_{x} c(x)$$
$$\text{subject to } \sum_{n=1}^{N} p_n \cdot \mathbb{I}(g(x, v_n) > 0) \leq p_f^t \tag{19}$$
$$x \in X,$$

where $\mathbb{I}(\cdot)$ denotes the indicator function. As discussed in Section 2.3, the computation of (19) is challenging because of the indicator function.

This issue can be addressed by replacing $p_f(x)$ with the buffered failure probability of the random variable $g(x, V)$, which we denote by $\bar{p}_f(x)$, and this produces the optimization problem

$$\min_{x} c(x)$$
$$\text{subject to } \bar{p}_f(x) \leq \bar{p}_f^t \tag{20}$$
$$x \in X.$$

The constraint in (20) can be reformulated in terms of a superquantile, i.e.,

$$\min_{x} c(x)$$
$$\text{subject to } \bar{q}_{\alpha_0}(x) \leq 0 \tag{21}$$
$$x \in X,$$

where $\alpha_0 = 1 - \bar{p}_f^t$. The replacement is valid because, given a solution $x$, $\bar{p}_f(x) \leq 1 - \alpha$ if and only if $\bar{q}_\alpha(x) \leq 0$ for $\alpha \in (0,1]$. Meanwhile, it is noted that given outcomes $g(x, v_n)$, $n = 1, \cdots, N$, a superquantile can be evaluated as an optimization problem (Rockafellar and Uryasev 2000):

$$\begin{aligned} \bar{q}_{\alpha_0}(x) &= \min_{z_0} z_0 + \frac{\mathbb{E}[\max\{g(x, V) - z_0, 0\}]}{1 - \alpha_0} \\ &= \min_{z_0} z_0 + \frac{1}{1 - \alpha_0} \sum_{n=1}^{N} p_n \cdot \max\{g(x, v_n) - z_0, 0\}. \end{aligned} \tag{22}$$

This expression can be incorporated into the optimization problem in (21) by introducing additional decision variables $z_0$ and $z = (z_1, \cdots, z_n)$ (Rockafellar and Royset 2010) as follows

$$\min_{x, z_0, z} c(x)$$

$$\text{subject to } z_0 + \frac{1}{\bar{p}_f^t} \sum_{n=1}^{N} p_n z_n \leq 0 \quad (23)$$

$$g(x, v_n) - z_0 \leq z_n, \; n = 1, \cdots, N$$

$$x \in X, \; z_0 \in \mathbb{R}, \; z_n \in \mathbb{R}^+, \; n = 1, \cdots, N,$$

where $\mathbb{R}$ and $\mathbb{R}^+$ denote the domain of real numbers and nonnegative real numbers, respectively. In the formulation, the first constraint controls the magnitude of the $(1 - \bar{p}_f^t)$-superquantile, while the constraints in the second line ensure that $z_n = \max\{g(x, v_n) - z_0, 0\}$, $n = 1, \cdots, N$, at optimality. In addition, the $\alpha_0$-quantile $q_{\alpha_0}$ of $g(x, V)$ is a minimizer of (22), so for an optimal $x^*$ in (23), one can identify the optimal $z_0$ in (23) as the $\alpha_0$-quantile of $g(x^*, V)$.

It is noted that (23), which is the final formulation that is used throughout the paper, is a mathematically well-behaved optimization problem in contrast to the formulation in (19) that includes the indicator function. The problem is linear with regard to $z_0$ and $z_n$, $n = 1, \cdots, N$, and therefore, the computational complexity only depends on the complexity of $c(x)$, $g(x, v_n)$, and $X$. For example, the problem becomes linear if $c(\cdot)$ is linear, $g(x, v_n)$ is linear in $x$, and $X$ is polyhedral. The problem becomes convex if $c(\cdot)$ is convex, $g(x, v_n)$ is convex in $x$, and $X$ is a convex set. Even when $c(x)$ and $g(x, v_n)$ are nonlinear or nonconvex functions, the optimization is still more efficient than the optimization with $p_f$ since the gradients of the constraints are well-defined as long as the limit-state function $g(x, V)$ has well-defined gradients with respect to $x$.

### 4.2  Active-set strategy for efficient optimization

In (23), the number of constraints linearly increases with the number of outcomes, $N$. In the case of sampling, $N$ is typically governed by the level of reliability and the desired coefficient of variance (c.o.v.) $\delta_t$ for the estimate of $p_f$ or $\bar{p}_f$. Specifically, the number of MCS samples required to meet $\delta_t$ is determined as

$$N = \frac{(1 - \bar{p}_f^t)}{\bar{p}_f^t \cdot \delta_t^2}. \quad (24)$$

In other words, $N$ increases with a lower $\bar{p}_f^t$ and $\delta_t$, producing large sample sizes in practice. This issue can be addressed by an *active-set* strategy (Example 6.16 of Royset and Wets (2021)), i.e., the optimization is performed only with a subset of the constraints in (23), namely the *active set*, that have actual influence on the optimal solution. This strategy is particularly effective when the number of constraints is much larger than that of the decision variables and, therefore, only a small subset of constraints are expected to have influence on the optimization result.

In the case of (23), the reliability constraint in the first line would typically be active, while the constraints in the second line are active only when the corresponding sample represents an event $g(x, v_n) - z_0 > 0$. Recall that $z_0$ eventually becomes equal to $q_{\alpha_0}$. Thus, this event corresponds indeed to those of concern for the buffered failure probability; see (9). Interestingly, the number of such failure events depends only on $\delta_t$ while being independent of $\bar{p}_f^t$ as

$$N_f = N \cdot \bar{p}_f^t = \frac{(1 - \bar{p}_f^t)}{\delta_t^2} \approx \frac{1}{\delta_t^2}, \quad (25)$$

where the final approximation is made from the common setting that $\bar{p}_f^t \ll 1$. As a result, the number of constraints can always be contained within a reasonable level, e.g., if $\delta_t = 0.05$, $N \cdot \bar{p}_f^t \approx 400$, which makes the optimization affordable even when using a general-purpose optimization solver. Since the set of failure events would change with varying decision variables $x$, the optimization needs to be carried out in an iterative fashion with updates of the active sets until the obtained solutions of $x$ converge.

The aforementioned strategy for BORM is summarized in **Algorithm 1.** The inputs of the algorithm are the initial solution $x_0$, the target parameters $\bar{p}_f^t$ and $\delta_t$, and a parameter $\beta$ that specifies the ratio of the number of active constraints to be considered relative to the number of failure events (e.g., 1.2). Thus, at a current solution $x^*$, the algorithm selects $N_a = \lceil \beta \cdot N_f \rceil$ samples with the largest values of $g(x^*, v_n)$ to constitute the constraints in the second line, where $N_f$ and $\lceil b \rceil$, respectively, refer to the number of the failure samples and the smallest integer that is greater than $b$. Then, by setting $z_n = 0$ for $n$ corresponding to the samples that are not selected, the optimization is performed to update the solution $x^*$. This process is iterated until the solution converges, e.g., $|x^* - x^{*\text{old}}|/|x^*| < \varepsilon$ with $x^{*\text{old}}$ and $\varepsilon$ respectively referring to the optimal solution of the previous iteration and a prespecified threshold.

It is noted that the active constraints of the final iteration can also be used to gain insights on the failure domain as they relate to the failure events. For example, they can be used to inspect the nonlinearity of the failure domain or to locate the most probable failure point (MPFP) similarly to the FORM and SORM.

---

**Algorithm 1.** BORM by active-set strategy

---

    **Procedure** BORM-Active-set (

        $x_0$      // initial solution

        $\bar{p}_f^t$, $\delta_t$      // target buffered failure probability and target c.o.v.

        $\beta$      // ratio of the number of active sets and the number of failure events

    )

1     $N \leftarrow \left\lceil \dfrac{(1 - \bar{p}_f^t)}{\bar{p}_f^t \cdot \delta_t^2} \right\rceil$

2     $N_f \leftarrow \lceil \bar{p}_f^t \cdot N \rceil$    // number of failure samples

3     $N_a \leftarrow \lceil \beta \cdot N_f \rceil$    // number of active samples

4     Randomly select data or generate MCS samples of size $N$, $v_n$, $n = 1, \cdots, N$

5     $x^* \leftarrow x_0$

6     While solution $x^*$ has not converged

7        $x^{*\text{old}} \leftarrow x^*$

8        Evaluate $g(x^*, v_n)$, $n = 1, \cdots, N$

9        Select $N_a$ samples with the largest values of $g(x^*, v_n)$

10       By including the selected samples and setting $z_n = 0$ for non-selected

           samples, optimize (23) with initial solution $x^{*\text{old}}$ to obtain new solution $x^*$

11    **return**  $x^*$

One of the practical issues in implementing the algorithm is that the solutions obtained at each iteration may show discontinuous leaps as the list of active constraints is updated. Such issue can be resolved by simple strategies that ensure the convergence of the solution. One of the possible approaches is to move toward the new solution from the old one by a limited step size. For example, in the numerical examples in Section 6, the solution $x^*$ is updated by interpolating between the new solution $x^{*\text{new}}$ and the old one $x^{*\text{old}}$ with a prespecified step size as

$$x^* \leftarrow x^{*\text{old}} + t^{h-1} \cdot \left(x^{*\text{new}} - x^{*\text{old}}\right) \qquad (26)$$

where $t$ and $h$ respectively denote the parameter of step size and the number of iterations. Although it is expected to be applicable in general, this strategy may not work for problems with discontinuous solution space, i.e., the interpolation between a pair of feasible solutions leads to infeasible solutions. In this case, the discontinuous leap of solutions can be discouraged by adding a penalty term to the cost function $c(x)$ as

$$c(x) + \lambda \|x - x^{*\text{old}}\|^2, \qquad (27)$$

where $\lambda > 0$ is a prespecified constant. However, it is noted that this strategy should be employed with caution as modifying the objective function might degrade the quality of the optimization results.

*4.3    Other issues for implementation*

Real-world problems often have multiple limit-state functions that are combined to form series, parallel, or general systems described by cut-sets and link-sets. A common approach is to assign a target probability to each limit-state function separately (Kim and Song 2021), i.e. $\bar{p}_{f,k}^t$ for each function $g_k(\cdot)$, $k = 1, \cdots, K$, which is also adopted for the numerical examples in Section 6. In this case, the optimization problem in (20) is expanded as

$$\min_{x} c(x)$$
$$\text{subject to } \bar{p}_{f,k}(x) \leq \bar{p}_{f,k}^t, \ k = 1, \cdots, K \qquad (28)$$
$$x \in X,$$

where $c(\cdot)$, $K$, and $X$ denote the cost function, the number of limit-state functions, and additional constraints imposed on the decision variables, respectively. In parallel to (21), this problem is equivalent to

$$\min_{x} c(x)$$
$$\text{subject to } \bar{q}_{\alpha_0,k}(x) \leq 0, \ k = 1, \cdots, K \qquad (29)$$
$$x \in X,$$

where $\bar{q}_{\alpha_0,k}(x)$ is the superquantile value associated with the $k$-th limit-state function. Thereby, the problem in (23) becomes expanded to

$$\min_{x, z_0, z} c(x)$$

$$\text{subject to } z_{0,k} + \frac{1}{\bar{p}_{f,k}^t} \sum_{n=1}^{N} p_n z_{n,k} \leq 0, k = 1, \cdots, K \qquad (30)$$

$$g_k(x, v_n) - z_{0,k} \leq z_{nk}, \ k = 1, \cdots, K \text{ and } n = 1, \cdots, N$$

$$x \in X, \ z_{0,k} \in \mathbb{R}, \ k = 1, \cdots, K, \text{ and } z_{nk} \in \mathbb{R}^+, \ k = 1, \cdots, K \text{ and } n = 1, \cdots, N,$$

where $z_0 = \{z_{0,1}, \cdots, z_{0,K}\}$; $z = \{z_{11}, \cdots, z_{1K}, \cdots, z_{N1}, \cdots, z_{NK}\}$. More detailed discussions, especially on other system types, can be found in Rockafellar and Royset (2010).

Most optimization tools benefit from being initialized with a good feasible solution. In the present case, a feasible solution can be found by relaxing the reliability constraint in (23) as

$$\min_{x, z_0, z} c(x) + \lambda \left( z_0 + \frac{1}{\bar{p}_f^t} \sum_{n=1}^{N} p_n z_n \right)$$

$$\text{subject to } g(x, v_n) - z_0 \leq z_n, \ n = 1, \cdots, N \qquad (31)$$

$$x \in X, \ z_0 \in \mathbb{R}, \ z_n \in \mathbb{R}^+, \ n = 1, \cdots, N,$$

where $\lambda > 0$ is a multiplier parameter; see Royset and Wets (2021), Section 5.D. If no feasible solution is found even with an arbitrarily large $\lambda$, it is concluded that the given problem is infeasible.

### 4.4 Refinements of BORM

Even with the active-set strategy, complicated formulas of the cost function and limit-state functions could result in a high computational cost of optimization. In this case, the optimization can be further facilitated by providing an initial solution that is near-optimal. Such solution can be found by performing a preliminary optimization with a higher target c.o.v. $\delta_{t,h}$ before computing the full-scale optimization with the original $\delta_t$. The optimization procedure is summarized as follows:

1. (If an initial solution is necessary but unknown) Obtain the initial solution by solving (31) with high target c.o.v. $\delta_{t,h}$ (e.g., 0.2)
2. Calibrate the initial solution by solving (23) with high target c.o.v. $\delta_{t,h}$
3. Using the calibrated initial solution, optimize (23) with the original target c.o.v. $\delta_t$ (e.g., 0.05)

## 5 Estimation of reliability sensitivity

Risk-informed decisions benefit from an assessment of which variables influence the failure probability most significantly. One of the common measures of such influence is the derivative of the failure probability with respect to the variable of interest, namely *reliability sensitivity* (Der Kiureghian et al., 2007). When using the conventional failure probability, estimating such derivative using discrete outcomes is challenging because of the indicator function as discussed in Section 2.3; see also Papaioannou et al. (2018). In contrast, the derivative of the buffered failure probability $\bar{p}_f$ does not involve the indicator function and can often be computed easily. Specifically, we compute the derivative in (11) with respect to any parameter $\theta$ defining $g(x, v_n)$, including components of the vectors $x$ and $v_n$, $n = 1, \cdots, N$. Suppose that $y_n = g(x, v_n)$, $n = 1, \cdots, N$, and $y_1 < y_2 < \cdots < y_N$. With $n^*$ satisfying (12) and $\partial y_n / \partial \theta = \partial g(x, v_n) / \partial \theta$, we obtain

$$\begin{aligned}\frac{\partial \bar{p}_f}{\partial \theta} &= \frac{\partial}{\partial \theta}\left(\frac{\sum_{n=n^*+1}^{N} p_n(y_n - y_{n^*})}{-y_{n^*}}\right) \\ &= \sum_{n=n^*+1}^{N} p_n \frac{\partial}{\partial \theta}\left(\sum_{n=n^*+1}^{N}\left(-\frac{y_n}{y_{n^*}} + 1\right)\right) \\ &= \sum_{n=n^*+1}^{N} p_n \left(\frac{\partial y_n}{\partial \theta}\left(-\frac{1}{y_{n^*}}\right) + \frac{\partial y_{n^*}}{\partial \theta}\left(\frac{y_n}{y_{n^*}^2}\right)\right).\end{aligned} \quad (32)$$

In the case of sampling (as discussed around (13)), an identical formula is available for the derivative of the estimate of $\hat{\bar{p}}_f$ of the buffered failure probability; then, $y_n$ would be computed using samples $\boldsymbol{v}_1, \cdots, \boldsymbol{v}_N$. We note that $\theta$ can be any parameter of the limit-state function as long as the derivative with respect to the parameter is well defined. This formula supplements those in Zhang et al. (2019), which apply under other hard-to-verify conditions. The present formula is especially useful due to its simplicity. In contrast, the derivative of the conventional failure probability is not informative under the present assumptions: $\partial p_f / \partial \theta = 0$ because small changes to $y_n$ does not change the failure probability.

## 6 Numerical examples

To demonstrate the efficiency and applicability of the proposed method, the optimization is computed using the general-purpose optimization solvers provided by Matlab®: fmincon in Section 6.1 and linprog and intlinprog in Section 6.2. While running fmincon, default settings are used except the maximum number of function evaluations and iterations, which are increased to $10^5$ and $10^4$, respectively, to prevent premature termination of the optimization. We use a personal desktop with processor Intel® Core™ i7 and RAM 16.0 GB.

*6.1 Benchmark examples*

This section addresses two benchmark examples that have been widely investigated in the literatures (Lee and Jung 2008; Kim and Song 2021). During the optimization, the uncertainties in the random variables are accounted for by generating MCS samples from the given distributions. Following the procedure proposed in Section 4.4, the target c.o.v.'s are set as $\delta_{t,h} = 0.2$ and $\delta_t = 0.05$ for the initial and full-scale computations, respectively, and the initial solutions are found using the penalized cost function in (31) with $\lambda = 10$. At each iteration, active sets are selected with $\beta = 1.2$ in **Algorithm 1**, and the solution has been updated using (26) with $t = 0.9$. For comparison with the literatures, the target buffered failure probability is set as the probability estimated from the optimal solutions provided in Kim and Song (2021). The probability is estimated using $10^6$ MCS samples, and since the examples have multiple limit-state functions, the maximum value is selected as the target value. In addition, after the optimization, the reliability sensitivity is estimated with respect to the parameters and decision variables. The accuracy of the estimation is confirmed as all results show difference less than 0.1 % from the numerical evaluations of $(\bar{p}_f(\theta + \Delta\theta) - \bar{p}_f(\theta))/\Delta\theta$ with $\Delta\theta = 10^{-6}$, where $\theta$ is the variable of interest.

6.1.1 A highly nonlinear limit-state function

This example has two decision variables $\boldsymbol{x} = (x_1, x_2)$ and the cost function

$$c(\boldsymbol{x}) = (x_1 - 3.7)^2 + (x_2 - 4)^2, \quad (33)$$

where $x_1 \in [0.00, 3.70]$ and $x_2 \in [0.00, 4.00]$, and these bounds specify the set $\mathbf{X}$. There are two random variables $\mathbf{V} = \{V_1, V_2\}$ and two limit-state functions (a highly nonlinear function and a linear one) that are defined as

$$g_1(\mathbf{x}, \mathbf{v}) = (x_1 + v_1)\sin(4(x_1 + v_1)) + 1.1(x_2 + v_2)\sin(2(x_2 + v_2)) \text{ and} \tag{34}$$
$$g_2(\mathbf{x}, \mathbf{v}) = -(x_1 + v_1) - (x_2 + v_2) + 3,$$

where $V_i \sim \mathcal{N}(0, \sigma_i^2)$, $i = 1,2$; $\sigma_1 = \sigma_2 = 0.1$; and $\mathcal{N}(\mu, \sigma^2)$ denotes the normal distribution with mean $\mu$ and standard deviation $\sigma$. Referring to the solution provided by Kim and Song (2021), the target buffered failure probability is set as $\bar{p}_f^t = 8.23 \times 10^{-2}$.

Table 1 summarizes the optimization results obtained from Kim and Song (2021), the solution using penalty parameter $\lambda = 10$, the preliminary solution with $\delta_{t,h}$, and the full-scale solution obtained with $\delta_t$. While the results generally agree with the results in Kim and Song (2021), it is noted that the final solution (obtained with $\delta_t$) leads to a higher $p_f$ than the reference solution, but a lower $\bar{p}_f$ and $\tau$. In other words, the use of $\bar{p}_f$ puts more emphasis on minimizing the heaviness of the upper tail than on minimizing the area beyond 0. This underlines the advantage of $\bar{p}_f$ for restraining high-consequence failure scenarios.

In Table 1, the efficiency of the proposed framework is evaluated by the number of iterations and the running time of **Algorithm 1**. The figures demonstrate the efficiency as the algorithm requires only a small number of iterations, especially for the full-scale optimization with $\delta_t$, and the total computational time is less than a minute. On the other hand, the nonlinearity of the problem is underlined in Figure 5, which illustrates the failure samples in the final round of optimization, i.e., $\mathbf{v}_n$ with $z_{n1}^* > 0$ or $z_{n2}^* > 0$ where $z_{n1}^*$ and $z_{n2}^*$, $n = 1, \cdots, N$, are the solutions from (30). It is noted that the proposed method enables a mathematically well-behaved optimization even when the limit-state is highly nonlinear, which is often challenging to address with conventional reliability methods.

From the optimal solution and the samples generated during the full-scale optimization, the reliability sensitivity is evaluated by the derivatives of $\bar{p}_f$ with respect to $\sigma_i$ and $x_i$, $i = 1,2$. The calculation is performed only for $g_1$ as $\bar{p}_f$ estimated for $g_2$ is zero. To evaluate (32), the derivatives of $g_1$ with respect to the variables are derived as

$$\frac{\partial g_1(\mathbf{x}, \mathbf{v}_n)}{\partial \sigma_1} = \frac{\partial g_1(\mathbf{x}, \mathbf{v}_n)}{\partial v_{1,n}} \cdot \frac{\partial v_{1,n}}{\partial \sigma_1}$$
$$= \left\{\sin\left(4(x_1 + v_{1,n})\right) + 4(x_1 + v_{1,n})\cos\left(4(x_1 + v_{1,n})\right)\right\} \cdot \frac{\partial v_{1,n}}{\partial \sigma_1},$$

$$\frac{\partial g_1(\mathbf{x}, \mathbf{v}_n)}{\partial \sigma_2} = \frac{\partial g_1(\mathbf{x}, \mathbf{v}_n)}{\partial v_2} \cdot \frac{\partial v_{2,n}}{\partial \sigma_2}$$
$$= \left\{1.1\sin\left(2(x_2 + v_{2,n})\right) + 2.2(x_2 + v_{2,n})\cdot \cos\left(2(x_2 + v_{2,n})\right)\right\} \cdot \frac{\partial v_{2,n}}{\partial \sigma_2}, \tag{35}$$

$$\frac{\partial g_1(\mathbf{x}, \mathbf{v}_n)}{\partial x_1} = \sin\left(4(x_1 + v_{1,n})\right) + 4(x_1 + v_{1,n})\cos\left(4(x_1 + v_{1,n})\right), \text{ and}$$

$$\frac{\partial g_1(\mathbf{x}, \mathbf{v}_n)}{\partial x_2} = 1.1\sin\left(2(x_2 + v_{2,n})\right) + 2.2(x_2 + v_{2,n})\cdot \cos\left(2(x_2 + v_{2,n})\right)$$

with

$$\frac{\partial v_{i,n}}{\partial \sigma_i} = \frac{\partial}{\partial \sigma_i}(\sigma_i u_{i,n}) = u_{i,n}, \quad i = 1,2, \tag{36}$$

where $u_{i,n}$ is the standard normal random variable that is used to generate sample $v_{i,n}$, $n = 1, \cdots, N$, which can be retrieved from $v_{i,n}$, i.e., $u_{i,n} = v_{i,n}/\sigma_i$. (Recall that $V_i \sim \mathcal{N}(0, \sigma_i^2)$). The results are summarized in Table 2, where it is found that $\bar{p}_f$ is more sensitive to $\sigma_1$ and $x_1$ than to $\sigma_2$ and $x_2$.

Table 1. Optimization results of the example with a highly nonlinear limit-state function: Solution in Kim and Song (2021), penalized cost function, preliminary optimization with c.o.v. 0.2, and full-scale optimization with c.o.v. 0.05.

| | | Kim and Song (2021) | Penalized cost function | Preliminary optimization | Full-scale optimization |
|---|---|---|---|---|---|
| Optimal solution, $(x_1, x_2)$ | | (2.81, 3.28) | (2.78, 2.52) | (2.87, 3.26) | (2.84, 3.26) |
| Cost | | 1.31 | 3.04 | 1.24 | 1.29 |
| Failure prob. estimate by $10^6$ samples, $(g_1, g_2)$ | $p_f$ | (0.0310, 0.00) | (0.00, 0.00) | (0.0528, 0.00) | (0.0327, 0.00) |
| | $\bar{p}_f$ | (0.0823, 0.00) | (0.00, 0.00) | (0.139, 0.00) | (0.0820, 0.00) |
| $\tau, (g_1, g_2)$ | | (2.65, N/A) | (N/A, N/A) | (2.63, N/A) | (2.51, N/A) |
| Number of iterations for active-set strategy | | - | 10 | 3 | 1 |
| Computation time (sec) | | - | 0.517 | 0.817 | 33.8 |

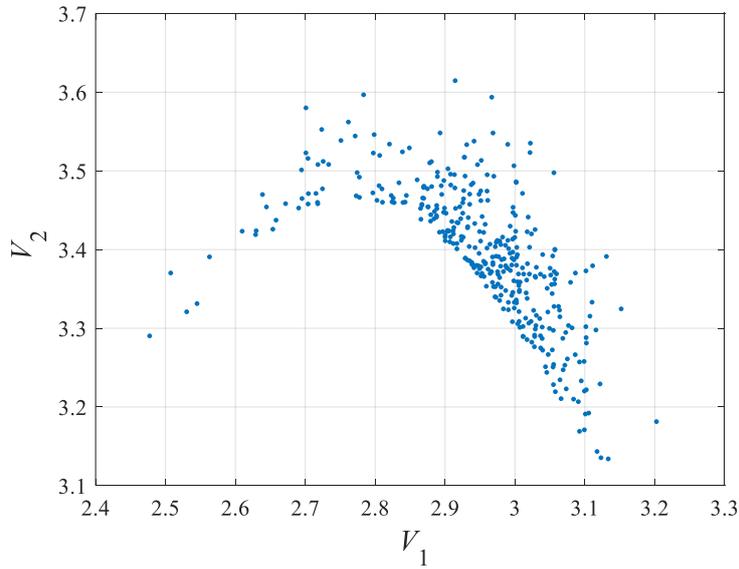

Figure 5. Samples of failure events in the example with a highly nonlinear limit-state function.

Table 2. Reliability sensitivity $\partial \bar{p}_f / \partial \theta$ of the example with a highly nonlinear limit-state function.

| | $\theta$ | $i = 1$ | $i = 2$ |
|---|---|---|---|
| $g_1$ | $\sigma_i$ | 2.31 | 1.60 |
| | $x_i$ | 1.32 | 1.23 |
| $g_2$ | $\sigma_i$ | N/A | N/A |
| | $x_i$ | N/A | N/A |

6.1.2 Welded beam structure

This example aims to find the optimal design of a welded beam structure which is illustrated in Figure 6 (Chen et al. 2013; Kim and Song 2021). There are four decision variables $x = (x_1, \cdots, x_4)$, and the cost function is defined as

$$c(x) = c_1 x_1^2 x_2 + c_2 x_3 x_4 (b_2 + x_2), \tag{37}$$

where $x_1 \in [3.175, 10]$, $x_2 \in [15, 254]$, $x_3 \in [200, 220]$, and $x_4 \in [3.175, 10]$; and these constraints define the set $X$. There are four random variables $V = \{V_1, \cdots, V_4\}$ and five limit-state functions given as

$$g_1(x, v) = \frac{\tau(x, v)}{b_6} - 1, \quad g_2(x, v) = \frac{\sigma(x, v)}{b_7} - 1, \quad g_3(x, v) = \frac{x_1 + v_1}{x_4 + v_4} - 1,$$

$$g_4(x, v) = \frac{\delta(x, v)}{b_5} - 1, \quad \text{and} \quad g_5(x, v) = 1 - \frac{P_c(x, v)}{b_1}, \tag{38}$$

where $V_i \sim \mathcal{N}(0, \sigma_i^2)$, $i = 1, \cdots, 4$, $\sigma_1 = \sigma_2 = 0.1693$, and $\sigma_3 = \sigma_4 = 0.0107$. The parameters of the cost and limit-state functions are summarized in Table 3, while the nested functions of the limit-state functions are given as

$$\tau(x, v) = \left[ t(x, v) + \frac{2t(x, v) tt(x, v)(x_1 + v_1)}{2R(x, v)} + tt(x, v)^2 \right]^{0.5},$$

$$t(x, v) = \frac{b_1}{\sqrt{2}(x_1 + v_1)(x_2 + v_2)}, \quad tt(x, v) = \frac{M(x, v) R(x, v)}{J(x, v)},$$

$$M(x, v) = b_1 (b_2 + 0.5(x_2 + v_2)),$$

$$R(x, v) = \sqrt{[(x_2 + v_2)^2 + \{(x_1 + v_1) + (x_3 + v_3)\}^2]/4}, \tag{39}$$

$$J(x, v) = \sqrt{2}(x_1 + v_1)(x_2 + v_2)[(x_2 + v_2)^2/12 + \{(x_1 + v_1) + (x_3 + v_3)\}^2/4],$$

$$\sigma(x, v) = \frac{6 b_1 b_2}{(x_3 + v_3)^2 (x_4 + v_4)}, \quad \delta(x, v) = \frac{4 b_1 b_2^3}{b_3 (x_3 + v_3)^3 (x_4 + v_4)}, \quad \text{and}$$

$$P_c(x, v) = \frac{4.013 (x_3 + v_3)(x_4 + v_4)^3 \sqrt{b_3 b_4}}{6 b_2^2} \left( 1 - \frac{x_3 + v_3}{4 b_2} \sqrt{\frac{b_3}{b_4}} \right).$$

From the optimal solution computed by Kim and Song (2021), the target buffered failure probability is set as $\bar{p}_f^t = 5.20 \times 10^{-3}$.

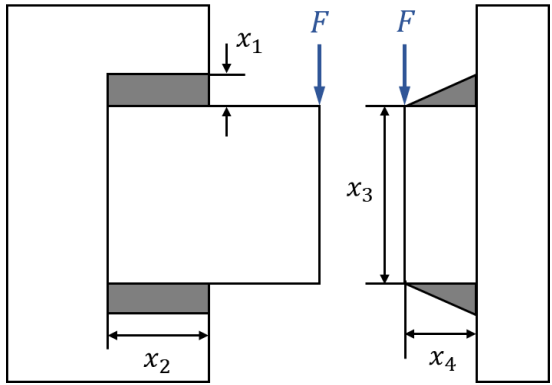

Figure 6. Welded beam structure (figure recreated from Chen et al. (2013)).

Table 3. Parameters of the cost function and limit-state functions of the welded beam example.

| Parameters | Value |
|---|---|
| $c_1$ | $6.74 \times 10^{-5}$ \$/mm³ |
| $c_2$ | $2.94 \times 10^{-6}$ \$/mm³ |
| $b_1$ | $2.67 \times 10^{4}$ N |
| $b_2$ | $3.56 \times 10^{2}$ mm |
| $b_3$ | $2.07 \times 10^{5}$ MPa |
| $b_4$ | $8.27 \times 10^{4}$ MPa |
| $b_5$ | 6.35 mm |
| $b_6$ | $9.38 \times 10$ MPa |
| $b_7$ | $2.07 \times 10^{2}$ MPa |

As reported in Table 4, the optimization results generally agree with the result in Kim and Song (2021). Nevertheless, there is a notable difference: Algorithm 1 produces slightly lower cost by better balancing the probabilities across the limit-state functions. In addition, as illustrated in the table, such mathematically better behaved properties do not compromise computational cost as the computation requires only 6 iterations and takes less than 5 minutes.

To estimate the reliability sensitivity, the derivatives $\partial \bar{p}_f / \partial \sigma_i$ and $\partial \bar{p}_f / \partial x_i$, $i = 1, \cdots, 4$, are evaluated using the obtained solution and samples of the final round. The derivatives of the limit-state functions with respect to the variables can be calculated in a similar way as in (35) and (36) of the previous example. The results are summarized in Table 5, which shows that $x_4$ and $s_4$ play the most significant role in determining the buffered failure probabilities of $g_2$, $g_3$, and $g_5$, while $x_1$ and $\sigma_1$ play the secondary role for $g_1$ and $g_3$.

Table 4. Optimization results of the welded beam example: Solution in Kim and Song (2021), penalized cost function, preliminary optimization with c.o.v. 0.2, and full-scale optimization with c.o.v. 0.05.

| | | Kim and Song (2021) | Penalized cost function | Preliminary optimization | Full-scale optimization |
|---|---|---|---|---|---|
| Optimal solution, $(x_1, \cdots, x_4)$ | | (5.72, 200, 211, 6.25) | (6.55, 242, 220, 10.0) | (5.76, 197, 211, 6.24) | (5.75, 199, 211, 6.24) |
| Cost | | 2.59 | 4.56 | 2.57 | 2.58 |
| Failure prob. estimate by $10^6$ samples, $(g_1, \cdots, g_5)$ | $p_f$ | (1.96, 0.00, 0.845, 0.00, 0.0110) × 10⁻³ | (0.00, 0.00, 0.00, 0.00, 0.00) | (2.86, 0.947, 2.29, 0.00, 0.984) × 10⁻³ | (1.95, 0.0920, 1.92, 0.00, 0.761) × 10⁻³ |
| | $\bar{p}_f$ | (5.20, 0.00, 2.25, 0.00, 0.0234) × 10⁻³ | (0.00, 0.00, 0.00, 0.00, 0.00) | (7.58, 2.47, 6.00, 0.00, 2.57) × 10⁻³ | (5.17, 0.229, 4.97, 0.00, 1.99) × 10⁻³ |
| $\tau$, $(g_1, \cdots, g_5)$ | | (2.65, N/A, 2.66, N/A, 2.12) | (N/A, N/A, N/A, N/A, N/A) | (2.65, 2.61, 2.61, N/A, 2.61) | (2.65, 2.49, 2.59, N/A, 2.61) |
| Number of iterations for active-set strategy | | - | 3 | 2 | 1 |
| Computation time (sec) | | - | 0.679 | 0.963 | 268 |

Table 5. Reliability sensitivity $\partial \bar{p}_f/\partial \theta$ of the welded beam example.

| | $\theta$ | $i = 1$ | $i = 2$ | $i = 3$ | $i = 4$ |
|---|---|---|---|---|---|
| $g_1$ | $\sigma_i$ | 0.285 | − 0.0000242 | − 0.0000573 | 0.00 |
| | $x_i$ | − 0.0975 | − 0.00186 | − 0.00215 | 0.00 |
| $g_2$ | $\sigma_i$ | 0.00 | 0.00 | − 0.00000685 | 0.413 |
| | $x_i$ | 0.00 | 0.00 | − 0.00672 | − 0.114 |
| $g_3$ | $\sigma_i$ | 0.265 | 0.00 | 0.00 | 0.0260 |
| | $x_i$ | 0.0930 | 0.00 | 0.00 | − 0.0930 |
| $g_4$ | $\sigma_i$ | N/A | N/A | N/A | N/A |
| | $x_i$ | N/A | N/A | N/A | N/A |
| $g_5$ | $\sigma_i$ | 0.00 | 0.00 | 0.000517 | 3.03 |
| | $x_i$ | 0.00 | 0.00 | − 0.00679 | − 0.997 |

## 6.2  Truss bridge system and weigh-in-motion (WIM) data

In order to demonstrate the applicability BORM for handling data-driven settings, this example finds the optimal cross section areas of the members in a truss bridge system against the traffic loads estimated by WIM data. As illustrated in Figure 7, the example truss bridge system has the length of 50 m and consists of 20 members. Considering the constructability, the members are classified into six groups as listed in Table 6; and within each group $i$, $i = 1,\cdots,6$, the cross section areas of the members are determined as $x_i$ such that $x_i \in [1.00, 4.00] \times 10^{-3}$ m² for all $i$. For optimization, the WIM data of 40,000 years are generated using the traffic loads model proposed by Kim and Song (2019). Then, using the influence lines, the dataset is translated into the samples $v_{mn}$ that represent the annual maximum internal force of member $m$ in year $n$, $m = 1,\cdots,20$ and $n = 1,\cdots,40{,}000$.

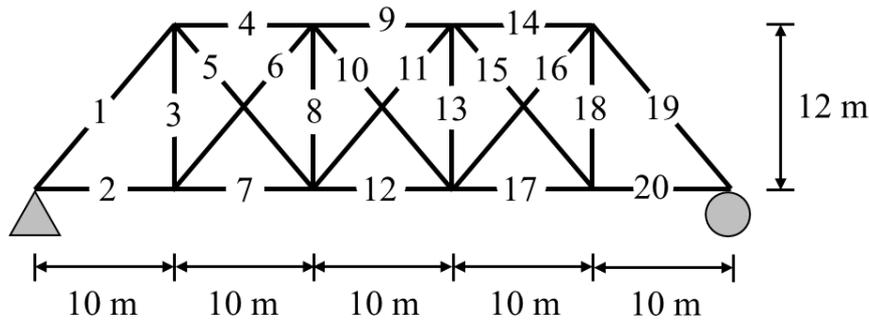

Figure 7. Example truss bridge system.

Table 6. Member types and the indices of associated members in the truss bridge example.

| Decision variable | Indices of associated members |
|---|---|
| $x_1$ | 1, 19 |
| $x_2$ | 2, 20 |
| $x_3$ | 3, 8, 13, 18 |
| $x_4$ | 4, 9, 14 |
| $x_5$ | 5, 6, 10, 11, 15, 16 |
| $x_6$ | 7, 12, 17 |

The cost function is set as the weight of the bridge, while the reliability constraints are defined such that, with respect to the self-weights and traffic loads, all members must have the probability of yielding less than 0.01/year, i.e., the return period of 100 years. Accordingly, the limit-state function of member $m$, $g_m(\pmb{x}, V_m)$ is formulated as

$$g_m(\pmb{x}, V_m) = \left| V_m + \sum_{i=1}^{6} x_i \cdot \left( \sum_{m' \in \mathcal{M}_i} l_{m'} w_{m'} \delta_{m'm} \right) \right| - E_{y,m} x_{i_m}, \tag{40}$$

where $\mathcal{M}_i$, $i = 1, \cdots, 6$, is the set of member indices associated with decision variable $x_i$; $i_m$ is the index such that $m \in \mathcal{M}_{i_m}$; and $l_m$ and $w_m = 7{,}950$ kg/m³, $m = 1, \cdots, 20$, are respectively the length and the unit weight of member $m$. Random variables $V_m$ and $E_{y,m}$ stand for the annual maximum internal force by traffic loads and the yield strength of member $m$. In the function, the first term represents the load demand which is the sum of traffic loads and self-weights, while the second term corresponds to the resistance force of the member. In the first term, $\delta_{m'm}$ refers to the internal force of member $m$ caused by the unit weight of member $m'$, which can be calculated by structural analysis.

The target buffered failure probability is evaluated from the proposed formula of $\tau$ in (16), i.e., $\bar{p}_f^t = p_f^t \times \tau^* = 0.01 \times 2.61$, leading to the optimization problem

$$\begin{aligned} \min_{\pmb{x}} c(\pmb{x}) &= \sum_{i=1}^{6} x_i \cdot \left( \sum_{m \in \mathcal{M}_i} l_m w_m \right) \\ \text{subject to } \bar{p}_{f,m}(\pmb{x}) &\leq 2.61 \times 10^{-2}, m = 1, \cdots, 20 \\ x_i &\in [1,4] \times 10^{-3}, i = 1, \cdots, 6. \end{aligned} \tag{41}$$

Then, using the samples of $V_m$ and $E_{y,m}$, $v_{mn}$ and $e_{y,mn}$, $m = 1, \cdots, 20$ and $n = 1, \cdots, N$, the optimization problem becomes

$$\begin{aligned} \min_{\pmb{x}, \pmb{z}_0, \pmb{z}, \pmb{s}} & \sum_{i=1}^{6} x_i \cdot \left( \sum_{m' \in \mathcal{M}_i} l_{m'} w_{m'} \right) \\ \text{subject to } & z_{0,m} + \frac{1}{\bar{p}_f^t} \frac{1}{N} \sum_{n=1}^{N} z_{mn} \leq 0, m = 1, \cdots, 20 \\ & h_{mn} - e_{y,mn} x_{i_m} - z_{0,m} \leq z_{mn}, \quad m = 1, \cdots, 6, \quad n = 1, \cdots, N \\ & v_{mn} + \sum_{i=1}^{6} x_i \cdot \left( \sum_{m' \in \mathcal{M}_i} l_{m'} w_{m'} \delta_{m'm} \right) \leq h_{mn}, \quad m = 1, \cdots, 6, n = 1, \cdots, N \\ & -\left\{ v_{mn} + \sum_{i=1}^{6} x_i \cdot \left( \sum_{m' \in \mathcal{M}_i} l_{m'} w_{m'} \delta_{m'm} \right) \right\} \leq h_{mn}, \quad m = 1, \cdots, 6, n = 1, \cdots, N \\ & x_i \in [1,4] \times 10^{-3}, i = 1, \cdots, 6, \ z_{0,m} \in \mathbb{R}, \ m = 1, \cdots, 20, \\ & z_{mn} \in \mathbb{R}^+, h_{mn} \in \mathbb{R}^+ \ m = 1, \cdots, 20, \ n = 1, \cdots, N. \end{aligned} \tag{42}$$

Interestingly, the problem is a linear programming (LP), which can be efficiently solved by general-purpose solvers.

Alternative, one might also consider the problem where the available cross section areas are limited to a finite set of choices. In this case, the problem becomes

$$\min_{x,z_0,z,s} \sum_{i=1}^{6} \left\{ \sum_{d=1}^{D_i} A_{id} x_{id} \cdot \left( \sum_{m' \in \mathcal{M}_i} l_{m'} w_{m'} \right) \right\}$$

$$\text{subject to } z_{0,m} + \frac{1}{\bar{p}_f^t} \frac{1}{N} \sum_{n=1}^{N} z_{mn} \leq 0, m = 1, \cdots, 20$$

$$h_{mn} - \sigma_{y,m} \sum_{d=1}^{D_{i_m}} A_{i_m d} x_{i_m d} - z_{0,m} \leq z_{mn}, \quad m = 1, \cdots, 20, \quad n = 1, \cdots, N$$

$$v_{mn} + \sum_{i=1}^{6} \left\{ \sum_{d=1}^{D_i} A_{id} x_{id} \cdot \left( \sum_{m' \in \mathcal{M}_i} l_{m'} w_{m'} \delta_{m'm} \right) \right\} \leq h_{mn}, \quad m = 1, \cdots, 20, n = 1, \cdots, N \quad (43)$$

$$-\left[ v_{mn} + \sum_{i=1}^{6} \left\{ \sum_{d=1}^{D_i} A_{id} x_{id} \cdot \left( \sum_{m' \in \mathcal{M}_i} l_{m'} w_{m'} \delta_{m'm} \right) \right\} \right] \leq h_{mn}, \quad m = 1, \cdots, 20, n = 1, \cdots, N$$

$$\sum_{d=1}^{D_i} x_{id} = 1, \ i = 1, \cdots, 6$$

$$x_{id} \in \{0,1\}, d = 1, \cdots, D_i, i = 1, \cdots, 6, \ z_{0,m} \in \mathbb{R}, \ m = 1, \cdots, 20,$$

$$z_{mn} \in \mathbb{R}^+, h_{mn} \in \mathbb{R}^+ \ m = 1, \cdots, 20, \ n = 1, \cdots, N,$$

where $D_i$ is the number of the candidate areas for member type $i$, $i = 1, \cdots, 6$, and $A_{id}$ is the area corresponding to solution $d$ of member type $i$, $d = 1, \cdots, D_i$. The decision variable $x_{id}$, $d = 1, \cdots, D_i$ and $i = 1, \cdots, 6$, is a binary variable indicating whether the solution $k$ is selected for member type $i$, while the last constraint ensures that a single solution is selected among the solutions $x_{i1}, \cdots, x_{iD_i}$. The problem in (43) is a mixed integer program (MIP), which can also be efficiently solved with general-purpose solvers.

The optimization is performed with three settings: (1) continuous decision variables $x$ and yield strength with no uncertainty, $\sigma_{y,m} \equiv 250$ MPa, for all $m = 1, \cdots, 20$, (2) continuous $x$ and $\sigma_{y,m} \sim \mathcal{N}(\mu_y, (\delta_y \mu_y)^2)$, for all $m$, where $\mu_y = 250$ MPa and $\delta_y = 0.1$, and (3) discrete $x$ and $\sigma_{y,m} \sim \mathcal{N}(\mu_y, (\delta_y \mu_y)^2)$, for all $m$. Table 7 summarizes the optimization results where it is found that the uncertainty in $\sigma_{y,m}$ and discrete $x$ lead to more conservative solutions, i.e., larger $x$. The computational cost of the optimization is also presented in the table, where all of the optimization takes less than 3 minutes. Figure 8 presents the probabilities $\bar{p}_f$ resulting from the obtained solutions, which are estimated by 40,000 samples. It is found that the failure is most likely to occur at members 1, 2, 3, 5, 9, and 12, for which the estimated values of $p_f$, $\bar{p}_f$, and $\tau$ under the second setting, i.e., continuous $x$ and uncertain $\sigma_{y,m}$, are summarized in Table 8. The buffered failure probabilities $\bar{p}_f$ are close to the target probability $\bar{p}_f^t = 2.61 \times 10^{-2}$ although the estimated values are slightly higher because of the numerical errors arising from sampling and optimization. The estimated values of $p_f$ are all close to 0.01, i.e., their return periods are around 100 years as desired. Meanwhile, the estimated $\tau$ values of members 2 and 3 are slightly higher than 2.72, which warns us that the members are more likely to have heavy-tailed distributions than other members.

The reliability sensitivity is estimated under the second setting, for which the derivatives of $g_m(x, V_m)$, $= 1, \cdots, 20$, with respect to $x_i$, $i = 1, \cdots, 6$, and $\delta_y$ are derived as

$$\frac{\partial g_m(\boldsymbol{x}, v_{mn})}{\partial x_i} = \begin{cases} h'_{mn} - \sigma_{y,mn}, & \text{if } m \in \mathcal{M}_i \\ h'_{mn}, & \text{otherwise} \end{cases} \quad (44)$$

where

$$h'_{mn} = \begin{cases} \sum_{m' \in \mathcal{M}_i} l_{m'} w_{m'} \delta_{m'm}, & \text{if } v_{mn} + \sum_{i=1}^{6} \left\{ \sum_{d=1}^{D_i} A_{id} x_{id} \cdot \left( \sum_{m' \in \mathcal{M}_i} l_{m'} w_{m'} \delta_{m'm} \right) \right\} \geq 0 \\ -\sum_{m' \in \mathcal{M}_i} l_{m'} w_{m'} \delta_{m'm}, & \text{otherwise} \end{cases} \quad (45)$$

and

$$\frac{\partial g_m(\boldsymbol{x}, v_{mn})}{\partial \delta_y} = \frac{\partial g_m(\boldsymbol{x}, v_{mn})}{\partial \sigma_{y,mn}} \cdot \frac{\partial \sigma_{y,mn}}{\partial \delta_y} = -x_{i_m} \cdot \frac{\partial \sigma_{y,mn}}{\partial \delta_y} \text{ and}$$
$$\frac{\partial \sigma_{y,mn}}{\partial \delta_y} = \frac{\partial}{\partial \delta_y}(\delta_y \mu_y u_{i,mn}) = \mu_y u_{i,mn}. \quad (46)$$

In the equations, $x_{i_m}$ is the decision variable with $m$ such that $m \in \mathcal{M}_{i_m}$; and $u_{i,mn}$ is the standard normal random variable that is used to generate sample $e_{y,mn}$, $m = 1, \cdots, 20$, $n = 1, \cdots, N$, which can be retrieved from $e_{y,mn}$ as $u_{i,n} = e_{y,mn}/(\delta_y \mu_y)$. Accordingly, the reliability sensitivity is estimated with respect to $\boldsymbol{x}$ and $\delta_y$ as illustrated in Figures 9(a) and (b), respectively. The results show that decision variable $x_i$ has the most significant effects on the members $m \in \mathcal{M}_i$ with the largest $\bar{p}_f$, while $\delta_y$ has greater influences on members with higher $\bar{p}_f$.

Table 7. Optimal solutions of the truss bridge example with (1) continuous decision variables, (2) continuous decision variables and uncertain yield strength, and (3) discrete decision variables and uncertain yield strength.

| Optimal solutions (x $10^{-3}$) | Continuous $\boldsymbol{x}$ | Continuous $\boldsymbol{x}$ and uncertain $\sigma_{y,m}$ | Discrete $\boldsymbol{x}$ and uncertain $\sigma_{y,m}$ |
|---|---|---|---|
| $x_1$ | 3.17 | 3.38 | 3.40 |
| $x_2$ | 2.03 | 2.17 | 2.80 |
| $x_3$ | 1.72 | 1.83 | 2.20 |
| $x_4$ | 3.06 | 3.26 | 3.40 |
| $x_5$ | 1.36 | 1.43 | 1.90 |
| $x_6$ | 2.83 | 3.06 | 3.10 |
| Computation time (sec) | 15.4 | 27.0 | 175 |

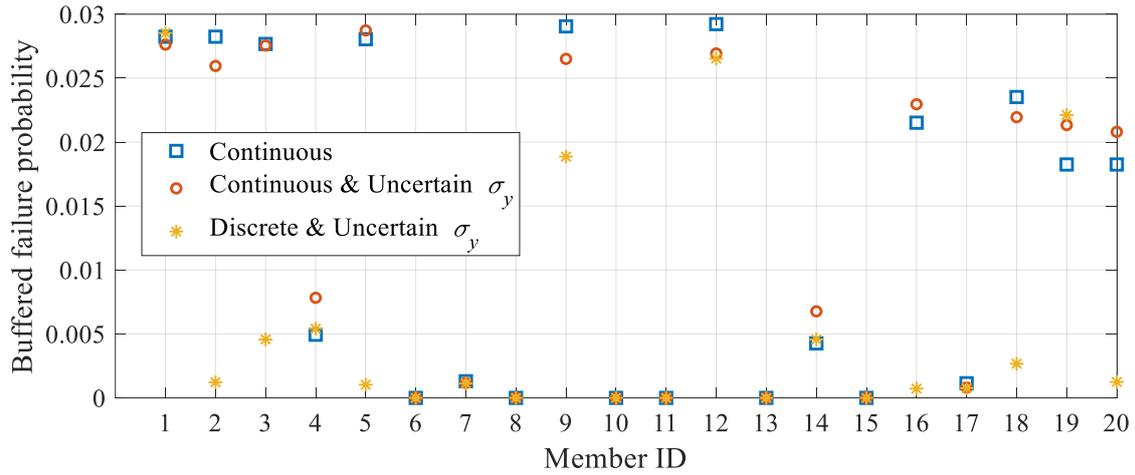

Figure 8. Buffered failure probabilities of the members of the truss bridge example with (1) continuous decision variables (blue squares), (2) continuous decision variables and uncertain yield strength (red circles), and (3) discrete decision variables and uncertain yield strength (yellow stars).

Table 8. Optimization results of the members with the highest failure probabilities, in the truss bridge example with continuous decision variables and uncertain yield strength.

| Continuous $x$ and uncertain $\sigma_{y,m}$ | Member index | | | | | |
|---|---|---|---|---|---|---|
| | 1 | 2 | 3 | 5 | 9 | 12 |
| $p_f$ (× 10⁻²) | 1.08 | 0.938 | 1.01 | 1.12 | 1.00 | 1.00 |
| $\bar{p}_f$ (× 10⁻²) | 2.76 | 2.60 | 2.76 | 2.87 | 2.65 | 2.69 |
| $\tau$ | 2.55 | 2.77 | 2.73 | 2.57 | 2.64 | 2.69 |

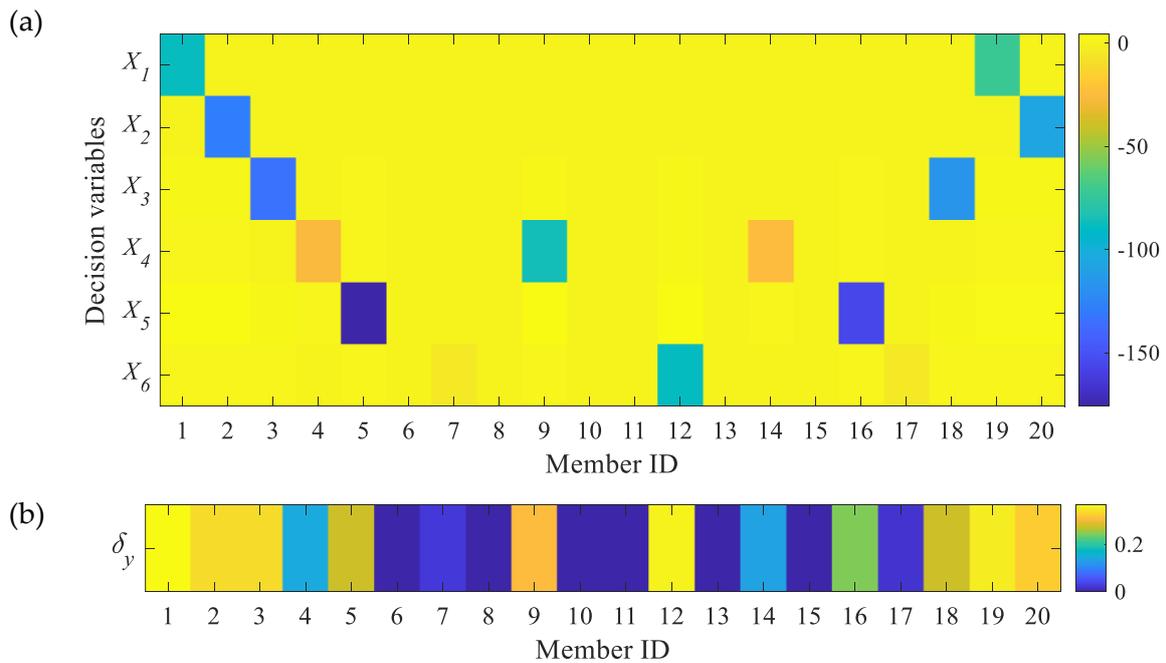

Figure 9. Reliability sensitivity of $\bar{p}_{f,m}$, $m = 1,\cdots,20$, with respect to (a) $x_i$, $i = 1,\cdots,6$, and (2) $\delta_y$ in the truss bridge example.

# 7 Conclusions

In order to take advantage of the recent advancement in data technology, this study proposes an efficient framework for data-driven optimization of reliability using the buffered failure probability, namely the *buffered optimization and reliability method* (BORM). The discussions demonstrate that the mathematically well-behaved properties of the buffered failure probability can improve the computational efficiency of inference tasks that otherwise would be challenging when using the conventional failure probability. To promote the applications of this alternative measure, which has rarely been employed for reliability engineering, this study presents novel expressions for the buffered failure probability and its sensitivities for settings with discrete distributions. Moreover, the characteristics of the buffered failure probability are systematically examined through common distributions, based on which an index is proposed as a measure of the heaviness of a distribution's tail, namely the *buffered tail index*. To facilitate the implementation of the index, reference values and expressions are proposed as well.

The proposed framework of BORM enables efficient data-driven reliability optimization. To this end, formulations and practical strategies of optimization are proposed, demonstrating the advantages of this alternative failure probability for reliability analysis using samples or data. The efficiency and applicability of the proposed optimization scheme are demonstrated by three numerical examples that involve highly nonlinear limit-state functions and high-dimensional distributions. All examples are solved using general-purpose optimization solvers and a personal desktop, which underlines the computational efficiency of the proposed optimization methodology. Specifically, the two benchmark examples demonstrate that compared to existing methods, BORM improves the mathematical properties of optimization problems and thereby leads to better solutions, while maintaining computational cost at practical levels; and the design optimization of a truss bridge demonstrates its utility for data-driven decision-making.

Using the buffered failure probability in parallel to the conventional failure probability has a great potential for solving the inference tasks that remain inefficient or even impossible by using the conventional probability alone, for which there are various topics that deserve further investigations. For example, while the numerical examples in this study only examine design optimization, other topics of optimization need to be explored as well. Such topics include sensor deployment, retrofitting, and system operations, for which data-driven decision-making would be even more relevant. Based on such developments, one can establish an automatic system of data collection and reliability optimization. Another issue worth investigating is the development of an efficient sampling technique for optimization using the buffered failure probability, so that the number of samples can be reduced, and thereby, the computation can become more efficient.

**Acknowledgment:** This work is supported in part by the Office of Naval Research under MIPR N0001421WX01496 and the Air Force Office of Scientific Research under MIPR F4FGA00350G004.

**Conflicts of Interest:** The authors declare no conflict of interest.

**Appendix A. Derivations of buffered tail index for common distributions**

With the buffered failure probability $\bar{p}_f = 1 - \alpha$, the conventional failure probability $p_f$ can be computed using $\alpha$-superquantile $\bar{q}_\alpha$ as $p_f = 1 - F(\bar{q}_\alpha)$ where $F(\cdot)$ refers to the cumulative distribution function (CDF). Accordingly, the buffered tail index $\tau = \bar{p}_f/p_f$ is derived as

$$\tau = \frac{\bar{p}_f}{p_f} = \frac{1-\alpha}{1-F(\bar{q}_\alpha)}. \tag{A1}$$

The closed-form expressions of $\alpha$-superquantile $\bar{q}_\alpha$ are derived for common distributions by Norton et al. (2019), which are summarized in Table A.1 for the normal, exponential, Weibull, and generalized extreme value (GEV) distributions. In the table, the probability density functions (PDFs) are presented as well to facilitate the understanding of the distribution parameters. Based on these formulas, the buffered tail index $\tau$ can be derived as illustrated in the table. It is noted that the index $\tau$ of the exponential distribution is a constant; $\tau$ of the normal distribution depends only on $\alpha$; and $\tau$ values of the other distributions depend on $\alpha$ as well as the parameters that determine the distribution shape, i.e., the standard deviation $s$ of the lognormal distribution, the shape parameter $k$ of the Weibull distribution, and the shape parameter $\xi$ of the GEV distribution.

Table A.1 PDF, $\alpha$-superquantile, and buffered tail index of common distributions

| Distribution | PDF $f(y)$, $\alpha$-superquantile $\bar{q}_\alpha$, and buffered tail index $\tau$ |
|---|---|
| Normal | $f(y) = \frac{1}{\sqrt{2\pi}} \exp\left(-\frac{(y-\mu)^2}{2\sigma}\right)$ where $\mu$ and $\sigma > 0$ are the mean and standard deviation. |
| | $\bar{q}_\alpha = \mu + \sigma \frac{\phi(\sqrt{2}\,\mathrm{erf}^{-1}(2\alpha-1))}{1-\alpha}$ where $\phi(\cdot)$ and $\mathrm{erf}^{-1}(\cdot)$ are respectively the PDF of the standard normal distribution and the inverse of the error function. |
| | $\tau = \dfrac{1-\alpha}{1-\Phi\left(\dfrac{\phi(\sqrt{2}\,\mathrm{erf}^{-1}(2\alpha-1))}{1-\alpha}\right)}$ where $\Phi(\cdot)$ is the CDF of the standard normal distribution. |
| Exponential | $f(y) = \begin{cases} \lambda \exp(-\lambda y), & y \geq 0 \\ 0, & \text{otherwise} \end{cases}$ where $\lambda > 0$ is the rate parameter. |
| | $\bar{q}_\alpha = \dfrac{-\ln(1-\alpha)+1}{\lambda}$ |
| | $\tau = e$ |
| Lognormal | $f(y) = \dfrac{1}{ys\sqrt{2\pi}} \exp\left(-\dfrac{(\ln y-\mu)^2}{2s^2}\right)$ where $\mu$ and $s > 0$ are the mean and standard deviation of $\ln y$. |
| | $\bar{q}_\alpha = \dfrac{1}{2} \exp\left(\mu + \dfrac{s^2}{2}\right) \cdot \dfrac{1+\mathrm{erf}\left(\dfrac{s}{\sqrt{2}} - \mathrm{erf}^{-1}(2\alpha-1)\right)}{1-\alpha}$ where $\mathrm{erf}(\cdot)$ is the error function. |
| | $\tau = \dfrac{1-\alpha}{1-\Phi\left(\left\{\dfrac{1}{s}\ln 2 + \dfrac{s}{2} + \dfrac{1}{s}\ln\left(\dfrac{1+\mathrm{erf}\left(\dfrac{s}{\sqrt{2}} - \mathrm{erf}^{-1}(2\alpha-1)\right)}{1-\alpha}\right)\right\}\right)}$ |

| | | |
|---|---|---|
| Weibull | $f(y) = \begin{cases} \dfrac{k}{\lambda}\left(\dfrac{y}{\lambda}\right)^{k-1} \exp\left(-\left(\dfrac{y}{\lambda}\right)^k\right), & y \geq 0 \\ 0, & \text{otherwise} \end{cases}$ where $k > 0$ and $\lambda > 0$ are the parameters of shape and scale, respectively. | |
| | $\bar{q}_\alpha = \dfrac{\lambda}{1-\alpha} \Gamma_U\left(1 + \dfrac{1}{k}, -\ln(1-\alpha)\right)$ where $\Gamma_U(a,b) = \int_b^\infty p^{a-1} \exp(-p)\, dp$ is the upper incomplete gamma function. | |
| | $\tau = \dfrac{1-\alpha}{\exp\left[-\left(\dfrac{1}{1-\alpha}\Gamma_U\left(1+\dfrac{1}{k},-\ln(1-\alpha)\right)\right)^k\right]}$ | |
| Generalized extreme value (GEV) | $f(y) = \begin{cases} \dfrac{1}{s}\left(1 + \dfrac{\xi(y-\mu)}{s}\right)^{-\frac{1}{\xi}-1} \exp\left[-\left(1 + \dfrac{\xi(y-\mu)}{s}\right)^{-\frac{1}{\xi}}\right], & \xi \neq 0 \\ \dfrac{1}{s}\exp\left(-\dfrac{y-\mu}{s}\right) \cdot \exp\left[-\exp\left(-\dfrac{y-\mu}{s}\right)\right], & \xi = 0 \end{cases}$ where $\mu$, $s > 0$, and $\xi$ are respectively the parameters of location, scale, and shape. | |
| | $\bar{q}_\alpha = \begin{cases} \mu + \dfrac{s}{\xi(1-\alpha)}\left[\Gamma_L\left(1-\xi, \ln\left(\dfrac{1}{\alpha}\right)\right) - (1-\alpha)\right], & \xi \neq 0 \\ \mu + \dfrac{s}{(1-\alpha)}(m + \alpha \ln(-\ln(\alpha)) - \text{li}(\alpha)), & \xi = 0 \end{cases}$ where $\Gamma_L(a,b) = \int_0^b p^{a-1}\exp(-p)\,dp$ is the lower incomplete gamma function; $\text{li}(x) = \int_0^x \dfrac{1}{\ln p}\,dp$ is the logarithmic integral function; and $m$ is the Euler-Mascheroni constant. | |
| | $\tau = \begin{cases} \dfrac{1-\alpha}{1-\exp\left[-\left(1+\dfrac{1}{(1-\alpha)}\left[\Gamma_L\left(1-\xi,\ln\left(\dfrac{1}{\alpha}\right)\right)-(1-\alpha)\right]\right)^{-1/\xi}\right]}, & \xi \neq 0 \\ \dfrac{1-\alpha}{1-\exp\left[-\exp\left\{-\dfrac{1}{(1-\alpha)}(m+\alpha\ln(-\ln(\alpha))-\text{li}(\alpha))\right\}\right]}, & \xi = 0 \end{cases}$ | |